\documentclass[12pt]{article}
\usepackage{fancyhdr,amsmath,amsfonts,amssymb}
\begin{document}
\setcounter{page}{23}
\vspace*{10mm}

\begin{center}
\uppercase{\textbf{Uniform model of geometric spaces}}
\end{center}

\bigskip

\begin{center}
\textsc{Alexandru Popa}
\end{center}

\bigskip

\textsc{Abstract.}
Full classification of geometric spaces was proposed by Isaak Yaglom in [2]. By defining the elliptic,
parabolic (or linear) and hyperbolic kinds of measure and applying them to distance, plane and dihedral
angle of different dimensions we get $3^n$ geometric spaces of dimension $n$. In his work [3] Yaglom
says that "finding a general description of all geometric systems [was] considered by mathematicians
the central question of the day." A. B. Khachaturean resumed Yaglom's work in [8].

Author developed a uniform model for all these spaces where distance and angle measure kinds are
parameters. This model is calculus centric, but can also be used in theoretical research. It is useful
in the following domains:
\begin{itemize}
\item deduction of uniform equations among geometric spaces;
\item uniform model applied to any space, which provides an easy way to calculate distances, plane and
dihedral angles of any dimension, areas and volumes as well as parallel (where applied) and orthogonal
property detection;
\item study of not yet described spaces and more.
\end{itemize}

\bigskip


\newpage
\lhead{}
\chead{Alexandru Popa - Uniform model of geometric spaces}
\rhead{}
\begin{center}
\textsc{1. Definitions}
\end{center}
As was shown by Yaglom in [2], some $n$-dimensional geometric space can be defined specifying its $n$
characteristics, or measure kinds. We will use numbers 1 for elliptic characteristic, 0 for parabolic
(or linear) one and $-1$ for hyperbolic one. So, full space specification is a set of $n$
characteristics $k_1,...,k_n \in\{-1, 0, 1\}$, which can be detected by a simple algorithm.

Define
\begin{equation}
K_i = \prod_{j=1}^ik_j,\,\forall i = \overline{0,n}.
\end{equation}
For two vectors $x, y \in \mathbb{PR}^n$, $x = \left<x_0:...:x_n\right>$, $y = \left<y_0:...:y_n\right>$
define a dot product in respect of characteristics $k_1...k_n$ as
\begin{equation}
x\odot y = \sum_{i=0}^nK_ix_iy_i.
\end{equation}
and cross product in respect of $k_1...k_n$ so that
\begin{equation*}
(x\odot y)^2 + k_1(x\otimes y)^2 = (x\odot x)(y\odot y),\, \forall x, y \in \mathbb{PR}^n.
\end{equation*}
It can be checked that\footnote{Here and further we will consider for simplicity that $k^0 = 1$ for
$k = 0$ too. We will say $x$ divide $k^i$, $k=0$ if in expression $x/k^i$ the exponent of $k$ in
numerator is greater then or equals to $i$.}
\begin{equation}
x\otimes y = \sqrt{\frac1{k_1}\sum_{i<j=0}^nK_iK_j(x_iy_j - x_jy_i)^2}.
\end{equation}
These products were considered by Klein in [1] for elliptic and hyperbolic spaces.

A $(n+1)\times(n+1)$ matrix is generalized orthogonal in respect of $k_1...k_n$ if for all columns
$c_i, c_j$ ($i, j = \overline{0,n}$)
\begin{equation}
\frac1{K_{min(i,j)}}c_i\odot c_j =
\begin{cases}
1, i = j,\\
0, i \ne j.
\end{cases}
\end{equation}

Having characteristics $k \in \{-1, 0, 1\}$ consider functions $C, S, T : \mathbb{R} \to \mathbb{R}$:
\begin{eqnarray}
C(x) = C(k, x) &=& \sum_{i=0}^\infty(-k)^i\frac{x^{2i}}{(2i)!},\\
S(x) = S(k, x) &=& \sum_{i=0}^\infty(-k)^i\frac{x^{2i+1}}{(2i+1)!},\\
T(x) = T(k, x) &=& \frac{S(k, x)}{C(k, x)}.
\end{eqnarray}
It's easy to see, that
\begin{eqnarray*}
C(x) =
\begin{cases}
\cos x, &k = 1,\\
1, &k = 0,\\
\cosh x, &k = -1
\end{cases}\\
S(s) =
\begin{cases}
\sin x, &k = 1,\\
x, &k = 0,\\
\sinh x, &k = -1
\end{cases}\\
T(x) =
\begin{cases}
\tan x, &k = 1,\\
x, &k = 0,\\
\tanh x, &k = -1
\end{cases}
\end{eqnarray*}

Define a geometric space with characteristics $k_1...k_n$ as "unit ball" in projective space:
$\mathbb{B}^n = \{x\in \mathbb{PR}^n\,|\, x\odot x=1\}$. Consider "points" $X\in\mathbb{B}^n$
corresponding vectors $x\in\mathbb{PR}^n$. Consider "space transformation" all linear mappings of
$\mathbb{PR}^n$ whose matrices are generalized orthogonal. They are also transformations of
$\mathbb{B}^n$ as they preserve it. Consider $m$-dimensional planes images of
$\mathbb{B}^m\subset \mathbb{B}^n$ on some transformation. All $m$-dimensional planes are (restricted
to $\mathbb{B}^n$) linear combination of first $m+1$ columns of some generalized orthogonal matrix.
So, we can identify $m$-dimensional planes, $m < n$ with such $(n+1)\times(m+1)$ matrices.

For two $m$-dimensional planes $X, Y$ define dot product in respect of $k_1...k_n$ as
\begin{equation}
X\odot Y = \sum_{i_0<...<i_m=0}^nX_{i_0...i_m}Y_{i_0...i_m}\prod_{p=1}^m\frac{K_{i_p}}{K_p},
\end{equation}
where
$$
M_{l_0...l_m} =
\begin{vmatrix}
m_{l_00}&\ldots&m_{l_0m}\\
\vdots&\ddots&\vdots\\
m_{l_m0}&\ldots&m_{l_mm}\\
\end{vmatrix}
$$
and cross product so that
\begin{equation*}
(X\odot Y)^2 + k_{m+1}(X\otimes Y)^2 = (X\odot X)(Y\odot Y)
\end{equation*}
It can be checked that
\begin{equation}
X\otimes Y = \sqrt{\frac1{k_{m+1}}\sum_{\substack{i_0<...<i_m=0\\j_0<...<j_m=0\\i_0...i_m<j_0...j_m}}^n
(X_{i_0...i_m}Y_{j_0...j_m} - X_{j_0...j_m}Y_{i_0...i_m})^2\prod_{p=1}^m\frac{K_{i_p}K_{j_p}}{K^2_p}}.
\end{equation}

This model generalizes spherical model of elliptic space, hyperboloid model of hyperbolic space [6],
projective euclidean space model [7] and describes many new spaces.

\begin{center}
\textsc{2. Calculus in uniform model}
\end{center}
Author shows that dot and cross products of points and planes is invariant in respect of space
transformation. Moreover, it can be used for distance and angle calculus based on equalities ($m < n$).
\begin{eqnarray}
X\odot Y &=& C_{m+1}(\phi),\\
X\otimes Y &=& S_{m+1}(\phi),
\end{eqnarray}
where $X$ and $Y$ are two points (if $m=0$) and $\phi$ is distance between them or $X$ and $Y$ are
$m$-dimensional planes (if $m>0$) and $\phi$ is angle between them and functions
$C_{m+1}(x) = C(k_{m+1}, x), S_{m+1}(x) = S(k_{m+1}, x)$.

For some figure $F \subset \mathbb{B}^n$ volume can be calculated using the following equation
\begin{equation}
V_{\mathbb{R}}(F) = \frac1{n+1}V_{\mathbb{B}}(C_F)
\end{equation}
where $C_F \subset \mathbb{R}^{n+1}$ is cone having origin $O = \{0,...,0\} \notin \mathbb{B}^n$ as
vertex and figure $F$ as base, $V_{\mathbb{B}}$ is native volume in $\mathbb{B}^n$ and $V_{\mathbb{R}}$
is volume in sense of $\mathbb{R}^{n+1}$. The advantage of this approach is the fact $V_{\mathbb{R}}$
is volume in a linear vector space which is usually easily to find.

Based on this unified model we can deduce common equation among all spaces. For example, consider
$\mathbb{B}^2$ with characteristics $k_1$ and $k_2$ and triangle $ABC\in\mathbb{B}^2$ with edges $a$,
$b$ and $c$, interior angles $\alpha$, $\gamma$ and exterior angle $\beta'$ (interior angle $\beta$ may
not exist). Then sine and cosine I and II lows have identical form in all 9 2-dimensional spaces:
\begin{equation}
\frac{S_1(a)}{S_2(\alpha)} = \frac{S_1(b)}{S_2(\beta')} = \frac{S_1(c)}{S_2(\gamma)},
\end{equation}
and
\begin{eqnarray}
C_1(a) &=& C_1(b)C_1(c) + k_1S_1(b)S_1(c)C_2(\alpha),\\
C_1(b) &=& C_1(a)C_1(c) - k_1S_1(a)S_1(c)C_2(\beta'),\\
C_1(c) &=& C_1(a)C_1(b) + k_1S_1(a)S_1(b)C_2(\gamma),
\end{eqnarray}
\begin{eqnarray}
C_2(\alpha) &=& C_2(\beta')C_2(\gamma) + k_2S_2(\beta')S_2(\gamma)C_1(a),\\
C_2(\beta') &=& C_2(\alpha)C_2(\gamma) - k_2S_2(\alpha)S_2(\gamma)C_1(b),\\
C_2(\gamma) &=& C_2(\alpha)C_2(\beta') + k_2S_2(\alpha)S_2(\beta')C_1(a),
\end{eqnarray}
or
\begin{eqnarray}
T_1^2(a) = \frac{T_1^2(b) + T_1^2(c) - 2T_1(b)T_1(c)C_2(\alpha) + k_1k_2T_1^2(b)T_1^2(c)S_1^2(\alpha)}
{(1 + k_1T_1(b)T_1(c)C_2(\alpha))^2},\\
T_1^2(b) = \frac{T_1^2(a) + T_1^2(c) + 2T_1(a)T_1(c)C_2(\beta') + k_1k_2T_1^2(a)T_1^2(c)S_1^2(\beta')}
{(1 - k_1T_1(a)T_1(c)C_2(\beta'))^2},\\
T_1^2(c) = \frac{T_1^2(a) + T_1^2(b) - 2T_1(a)T_1(b)C_2(\gamma) + k_1k_2T_1^2(a)T_1^2(b)S_1^2(\gamma)}
{(1 + k_1T_1(a)T_1(b)C_2(\gamma))^2},
\end{eqnarray}
\begin{eqnarray}
T_2^2(\alpha) = \frac{T_2^2(\beta') + T_2^2(\gamma) - 2T_2(\beta')T_2(\gamma)C_1(a)
+ k_1k_2T_2^2(\beta')T_2^2(\gamma)S_1^2(a)}{(1 + k_2T_2(\beta')T_2(\gamma)C_1(a))^2},\\
T_2^2(\beta') = \frac{T_2^2(\alpha) + T_2^2(\gamma) + 2T_2(\alpha)T_2(\gamma)C_1(b)
+ k_1k_2T_2^2(\alpha)T_2^2(\gamma)S_1^2(b)}{(1 - k_2T_2(\alpha)T_2(\gamma)C_1(b))^2},\\
T_2^2(\gamma) = \frac{T_2^2(\alpha) + T_2^2(\beta') - 2T_2(\alpha)T_2(\beta')C_1(c)
+ k_1k_2T_2^2(\alpha)T_2^2(\beta')S_1^2(c)}{(1 + k_2T_2(\alpha)T_2(\beta')C_1(c))^2}.
\end{eqnarray}

As another example, consider $\mathbb{B}^2$ with characteristics $k_1, k_2 = 1$ and $ABC\in\mathbb{B}^2$
right triangle with catheti $a, b$, hypotenuse $c$ and angles $\alpha$ and $\beta$. Equations of $ABC$
have the same form for elliptic, euclidean and hyperbolic planes.
\begin{eqnarray}
T_1^2(c) &=& T^2_1(a) + T_1^2(b) + k_1T^2_1(a)T_1^2(b),\\
T_1(b) &=& T_1(c)\cos\alpha,\\
T_1(a) &=& T_1(c)\cos\beta,\\
S_1(a) &=& S_1(c)\sin\alpha,\\
S_1(b) &=& S_1(c)\sin\beta,\\
T_1(a) &=& S_1(b)\tan\alpha,\\
T_1(b) &=& S_1(a)\tan\beta,\\
\cos\alpha &=& C_1(a)\sin\beta,\\
\cos\beta &=& C_1(b)\sin\alpha,\\
C_1(c) &=& \cot\alpha\cot\beta.
\end{eqnarray}

\bigskip
\begin{center}
\textsc{References}
\end{center}

[1] Felix Klein, Vorlesungen Nicht-Euklidische Geometrie, B.G.Teubner, Leipzig 1890.

[2] Isaak Yaglom, A simple non-euclidean geometry and its physical basis, Springer, New York 1979.

[3] Isaak Yaglom, Felix Klein and Sophus Lie, Birkhauser, 1988.

[4] Fenchel, Werner, Elementary geometry in hyperbolic space, De Gruyter Studies in mathematics. 11.
Berlin-New York: Walter de Gruyter \& Co 1989.

[5] Naber, Gregory L., The Geometry of Minkowski Spacetime. New York, Springer-Verlag 1992,
ISBN 0387978488.

[6] Reynolds, William F, Hyperbolic Geometry on a Hyperboloid, American Mathematical Monthly 1993,
100:442-455.

[7] Coxeter H. S. M., The Real Projective Plane, 3rd ed, Springer Verlag 1995.

[8] A. B. Khachaturean, Galilean geometry, MCNMO, Moskow, 2005.

[9] James W. Anderson, Hyperbolic Geometry, Springer 2005, ISBN 1852339349

\bigskip

\noindent 
Alexandru Popa\\
Department of Computer Sciences\\
Vest University of Timisoara\\
Address: Blvd. V. Parvan 4, Timisoara 300223, Timis, Romania\\
email:\textit{alpopa@gmail.com}
\end{document}